# A New Approach in Plane Kinematics

*Symplectic Kinematics in $\mathbb{R}^2$*

Stefan Gössner[1]

[1]Dortmund University of Applied Sciences. Department of Mechanical Engineering

## *Abstract*

The kinematics of particles and rigid bodies in the plane are investigated up to higher-order accelerations. Discussion of point trajectories leads from higher-order poles to higher-order Bresse circles of the moving plane. Symplectic geometry in vector space $\mathbb{R}^2$ is used here as a new approach and leads to some new recursive vector formulas. This article is dedicated to the memory of Professor Pennestri.

# Content

# 1. Introduction and Prerequisites

Kinematics is Geometry in Motion. A distinction is generally made between the kinematics of particles and the kinematics of rigid bodies [1,3]. In technical textbooks, motion is often categorised into 2D motion and 3D motion [1,3,4,15]. 3D motion is usually described using vector algebra with transformation matrices [1,3,4,15].

Plane kinematics continues to be of great importance in engineering. It turns out that it cannot be treated as a special case of 3D vector algebra by simply setting the third coordinate to zero. Kinematics in the plane is sometimes described using complex numbers instead of vector algebra [4,5,6,15,17]. It appears that the vector space $\mathbb{R}^2$ gets isomorphic to complex numbers by adding a complex structure. This leads us to *symplectic geometry* [1,8,16], which is largely coordinate-, trigonometry- and matrix-free and offers advantages over complex numbers in planar geometric applications.

This paper is structured as follows. We begin with the kinematics of particles in the plane, describing their position and higher-order time derivatives using a single vector in each case. With this knowledge, we take a closer look at the trajectory of a particle and its properties. In the following section, we will discuss rigid body kinematics, which will lead us further to the higher-order poles of plane motion. After a digression into relative kinematics, higher-order Bresse circles are investigated. The last section deals exclusively with geometric aspects of planar motion. Here, the Euler-Savary equation, polodes and their higher derivatives, as well as their curvature properties, are examined in more detail. The discussion of Bottema invariants concludes this section.

The analyses carried out lead to generally known results. However, some of the recursive vector formulas are new to the best of the author's knowledge and belief.

## 1.1 Symplectic Geometry in $\mathbb{R}^2$

Starting with the *Euclidean vector space* we get the *standard scalar product* (Euclidean *structure g*), which associates a number to every pair of vectors $\begin{pmatrix} a_1 \\ a_2 \end{pmatrix}$ and $\begin{pmatrix} b_1 \\ b_2 \end{pmatrix}$

$$g(\mathbf{a}, \mathbf{b}) = \mathbf{a}\mathbf{b} = a_1 b_1 + a_2 b_2 \,.$$

Then we are adding a *complex structure* $\mathbf{J} = \begin{pmatrix} 0 & -1 \\ 1 & 0 \end{pmatrix}$ with $\mathbf{J}^{-1} = \mathbf{J}^T = -\mathbf{J}$ and $\mathbf{J}^2 = -\mathbf{I}$. Complex structure $\mathbf{J}$ is an *orthogonal operator* and transforms any vector into a *skew-orthogonal* one [1,8,16]. Thus we entered the *complex vector space*. As a shortcut we will place a *tilde* '~' symbol over the skew-orthogonal vector variable.

$$\tilde{\mathbf{a}} = \mathbf{J}\mathbf{a} = \begin{pmatrix} -a_2 \\ a_1 \end{pmatrix} \,.$$

Whenever you encounter a vector symbol with a tilde $\tilde{\mathbf{x}}$ in the following, you can replace it with $\mathbf{J}\mathbf{x}$. Yet applying the orthogonal operator to the first vector in the scalar product above gets us to the *skew-scalar product* (*symplectic structure $\omega$*)

$$\omega(\mathbf{a}, \mathbf{b}) = \tilde{\mathbf{a}}\mathbf{b} = a_1 b_2 - a_2 b_1 \,.$$

The skew-scalar product gives us the *area* of the parallelogram spanned from vector $\mathbf{a}$ to vector $\mathbf{b}$, which is a *directed* or *oriented* area [1,8].

When using the orthogonal operator '~', following rules apply:

$$\widetilde{\mathbf{a} + \mathbf{b}} = \tilde{\mathbf{a}} + \tilde{\mathbf{b}}\,; \quad \tilde{\tilde{\mathbf{a}}} = -\mathbf{a}\,; \quad \tilde{\mathbf{a}}\mathbf{a} = 0\,; \quad \tilde{\mathbf{a}}\mathbf{b} = -\mathbf{a}\tilde{\mathbf{b}}\,; \quad \tilde{\mathbf{a}}\tilde{\mathbf{b}} = \mathbf{a}\mathbf{b}$$

Since the skew-scalar product from $\mathbf{b}$ to $\tilde{\mathbf{a}}$ and from $\mathbf{a}$ to $\tilde{\mathbf{b}}$ is equivalent to the standard scalar product of $\mathbf{a}$ and $\mathbf{b}$ according to the last rule, i.e. $\tilde{\mathbf{b}}\tilde{\mathbf{a}} = \tilde{\mathbf{a}}\tilde{\mathbf{b}} = \mathbf{a}\mathbf{b}$, we can interpret the standard scalar product as a parallelogram area as well (Fig. 1a).

> *Symplectic geometry is an areal type of geometry. So areas are first class citizens instead of lengths and angles from Euclidean geometry* [8].
> (Dusa Mc. Duff)

The Euclidean, complex, and symplectic structures together are referred to as a *compatible triple*. If two of these are given, the third is automatically defined. Ultimately, we are working with three compatible, harmonizing vector spaces – the Euclidean, complex, and symplectic vector space in $\mathbb{R}^2$ [16].

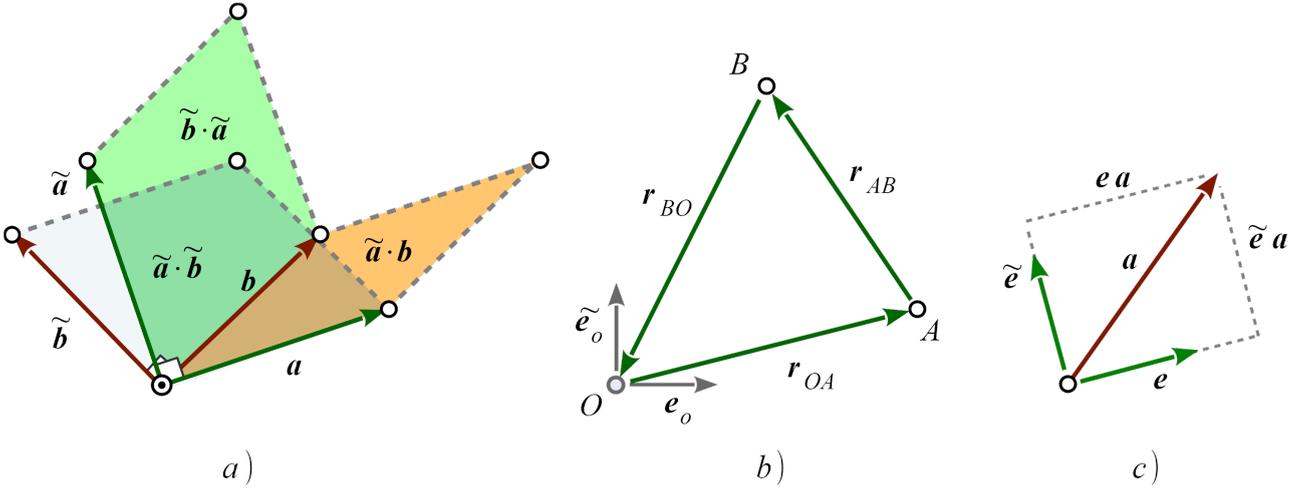

Fig. 1: Characteristics of symplectic vectors.

Symplectic geometry is fundamentally coordinate-free. One can interpret the components $\begin{pmatrix} a_1 \\ a_2 \end{pmatrix}$ of vector $\mathbf{a}$ with respect to an implicit global orthonormal symplectic basis $\{\mathbf{e}_o, \tilde{\mathbf{e}}_o\}$, where $\mathbf{e}_o$ is directed horizontally to the right (Fig. 1b). In the same way we can reformulate the components of vector $\mathbf{a}$ with respect to any arbitrary symplectic basis $\{\mathbf{e}, \tilde{\mathbf{e}}\}$ by $\begin{pmatrix} \mathbf{ea} \\ \tilde{\mathbf{e}}\mathbf{a} \end{pmatrix}$ (Fig. 1c).

Vectors are free. They are not bound to a fixed global coordinate origin. Points are denoted by capital letters in this article. They exist conceptually, but are not vectors here, unlike points as complex numbers. However, the vector from one point to another is defined.

> **Definition 1.1:**
> Given two points $A$ and $B$, the vector $\mathbf{r}_{AB}$ describes the transition from $A$ to $B$. Swapping the indices results in $\mathbf{r}_{BA} = -\mathbf{r}_{AB}$.

An explicit origin $O$ is not needed. If we want to, we are free to define any point on the plane as the origin $O$ (Fig. 1b). The vector closure equation of the triangle is then $\mathbf{r}_{OA} + \mathbf{r}_{AB} + \mathbf{r}_{BO} = \mathbf{0}$.

As an exception to the rule of always writing two point indices, we are allowed to omit the index $O$ of the origin if and only if it is in the first position. The above closure equation can therefore be rewritten as $\mathbf{r}_{AB} = \mathbf{r}_{OB} - \mathbf{r}_{OA} = \mathbf{r}_B - \mathbf{r}_A$. Note, that the irrelevance of points as bound vectors contrasts with the relevance of their time derivatives (velocity, acceleration, ...).

## 1.2 Rotation in the Plane

The rotation matrix $\mathbf{R}(\theta)$ is an *orthonormal matrix* with $\mathbf{R}^T = \mathbf{R}^{-1}$, which can be decomposed into a symmetric and a skew-symmetric part

$$\mathbf{R}(\theta) = \begin{pmatrix} \cos\theta & -\sin\theta \\ \sin\theta & \cos\theta \end{pmatrix} = \cos\theta \cdot I + \sin\theta \cdot \mathbf{J}\,.$$

This gives us the symplectic equivalence to Euler's formula $e^{ix} = \cos x + i \sin x$ in $\mathbb{R}^2$. So rotating a vector $\mathbf{a}$ by angle $\theta$ can be formulated in a matrix-free notation

$$\mathbf{Ra} = \cos\theta\, \mathbf{a} + \sin\theta\, \tilde{\mathbf{a}} \qquad (1.1)$$

where the expression on the right side conforms to the complex product.

> **Lemma 1.2:**
> The $k$-order derivative of the rotation matrix $\mathbf{R}$ with respect to angle $\theta$ is
>
> $$\mathbf{R}^{(k)} = \frac{d^k \mathbf{R}}{d\theta^k} = \mathbf{R}\mathbf{J}^k \quad \text{with} \quad \mathbf{J}^k = \begin{cases} \mathbf{I} & \text{for } k = 0 \pmod 4 \\ \mathbf{J} & \text{for } k = 1 \pmod 4 \\ -\mathbf{I} & \text{for } k = 2 \pmod 4 \\ -\mathbf{J} & \text{for } k = 3 \pmod 4 \end{cases} \quad (1.2)$$

*Proof.* Deriving equation (1.1) with respect to $\theta$ (denoted by prime) yields

$$\mathbf{R}'\mathbf{a} = -\sin\theta\,\mathbf{a} + \cos\theta\,\tilde{\mathbf{a}} = \mathbf{J}\mathbf{R}\mathbf{a}$$

with the insight that $\mathbf{R}' = \mathbf{R}\mathbf{J} = \mathbf{J}\mathbf{R}$ applies. Via $\mathbf{R}'' = \mathbf{R}'\mathbf{J} = \mathbf{R}\mathbf{J}^2$ and continuing, we finally arrive at relation (1.2). □

> **Proposition 1.3:**
> Rotating two vectors by the same angle has no effect on their (skew-)scalar product, i.e. $(\mathbf{R}\mathbf{a})(\mathbf{R}\mathbf{b}) = \mathbf{a}\mathbf{b}$.

## 1.3 Kinematics in the Plane

Planar kinematics distinguishes between bodies without extension and those with extension. In the first case, we consider the motion of a particle with two degrees of freedom. Its position, velocity and acceleration are each described by a single vector.

In the second case of a rigid body, two points or one point plus a direction are sufficient to uniquely describe its motion with three degrees of freedom. The contour of the body is mostly irrelevant. That is why we often refer to *planes* without boundaries moving relative to other (parallel) planes. Each plane carries an implicit frame – a symplectic basis – for quantifying relative motion parameters of that plane. These motion parameters are called *invariants* because they do not depend on specific points on that plane.

# 2. Kinematics of a Single Vector in $\mathbb{R}^2$

It is advantageous to discuss the motion of points using polar vectors, as these are easy to differentiate. The polar components of a single vector $\mathbf{r}$ are its strictly separated values, magnitude $r$ and unit vector $\mathbf{e}$.

$$\mathbf{r} = r\,\mathbf{e} = r\begin{pmatrix} \cos\varphi \\ \sin\varphi \end{pmatrix} \quad (2.1)$$

Deriving unit vector $\mathbf{e}$ with respect to $\varphi$ gives

$$\mathbf{e}' = \begin{pmatrix} -\sin\varphi \\ \cos\varphi \end{pmatrix} = \mathbf{J}\mathbf{e} = \tilde{\mathbf{e}} \quad \text{and} \quad \mathbf{e}^{(k)} = \frac{d^k \mathbf{e}}{d\varphi^k} = \mathbf{J}^k \mathbf{e},$$

with $\mathbf{J}^k$ according to (1.2).

> **Theorem 2.1**:
> The $k+1$ time derivative of a vector in the plane is a linear recurrence relation of its previous *parallel* and *orthogonal* $k$th time derivative components with respect to its symplectic basis $\{\mathbf{e}, \tilde{\mathbf{e}}\}$
>
> $$\mathbf{r}^{(k+1)} = \underbrace{(\dot{a}_{\parallel}^{(k)} - \dot{\varphi} a_{\perp}^{(k)})}_{a_{\parallel}^{(k+1)}} \mathbf{e} + \underbrace{(\dot{a}_{\perp}^{(k)} + \dot{\varphi} a_{\parallel}^{(k)})}_{a_{\perp}^{(k+1)}} \tilde{\mathbf{e}} \qquad (2.2)$$

*Proof.* The first two time derivatives of equation (2.1) are

$$\dot{\mathbf{r}} = \dot{r}\,\mathbf{e} + \dot{\varphi}\,r\,\tilde{\mathbf{e}}$$
$$\ddot{\mathbf{r}} = (\ddot{r} - \dot{\varphi}^2 r)\,\mathbf{e} + (2\dot{\varphi}\dot{r} + \ddot{\varphi} r)\,\tilde{\mathbf{e}}$$

For the general case of the derivative $\mathbf{r}^{(k)}$ we formulate it as a linear combination with respect to the symplectic basis $\{\mathbf{e}, \tilde{\mathbf{e}}\}$.

$$\mathbf{r}^{(k)} = a_{\parallel}^{(k)} \mathbf{e} + a_{\perp}^{(k)} \tilde{\mathbf{e}}$$

Differentiating that with respect to time gives us

$$\mathbf{r}^{(k+1)} = \dot{a}_{\parallel}^{(k)} \mathbf{e} + a_{\parallel}^{(k)} \dot{\varphi} \tilde{\mathbf{e}} + \dot{a}_{\perp}^{(k)} \tilde{\mathbf{e}} - a_{\perp}^{(k)} \dot{\varphi} \mathbf{e}$$

Summerizing yields equation (2.2). $\square$

Calculating the first five time derivatives of $a_{\parallel}^{(k)}$ and $a_{\perp}^{(k)}$ by recursion (2.2) leads to

$$
\begin{aligned}
a_{\parallel}^{(0)} &= r\,; & a_{\perp}^{(0)} &= 0 \\
a_{\parallel}^{(1)} &= \dot{r}\,; & a_{\perp}^{(1)} &= \dot{\varphi} r \\
a_{\parallel}^{(2)} &= \ddot{r} - \dot{\varphi}^2 r\,; & a_{\perp}^{(2)} &= 2\dot{\varphi}\dot{r} + \ddot{\varphi} r \\
a_{\parallel}^{(3)} &= \dddot{r} - 3\dot{\varphi}^2 \dot{r} - 3\ddot{\varphi}\dot{\varphi} r\,; & a_{\perp}^{(3)} &= 3\dot{\varphi}\ddot{r} + 3\ddot{\varphi}\dot{r} + (\dddot{\varphi} - \dot{\varphi}^3) r \\
a_{\parallel}^{(4)} &= \ddddot{r} - 6\dot{\varphi}^2 \ddot{r} - 12\ddot{\varphi}\dot{\varphi}\dot{r} - (4\dddot{\varphi}\dot{\varphi} + 3\ddot{\varphi}^2 - \dot{\varphi}^4) r\,; & a_{\perp}^{(4)} &= 4\dot{\varphi}\dddot{r} + 6\ddot{\varphi}\ddot{r} + 4(\dddot{\varphi} - \dot{\varphi}^3)\dot{r} + (\ddddot{\varphi} - 6\ddot{\varphi}\dot{\varphi}^2) r \\
a_{\parallel}^{(5)} &= \dddddot{r} - 10\dot{\varphi}^2 \dddot{r} - 30\ddot{\varphi}\dot{\varphi}\ddot{r} - (20\dddot{\varphi}\dot{\varphi} + 15\ddot{\varphi}^2 - 5\dot{\varphi}^4)\dot{r}\,; & a_{\perp}^{(5)} &= 5\dot{\varphi}\ddddot{r} + 10\ddot{\varphi}\dddot{r} + 10(\dddot{\varphi} - \dot{\varphi}^3)\ddot{r} + (5\ddddot{\varphi} - 30\ddot{\varphi}\dot{\varphi}^2)\dot{r} \\
&\quad -(5\ddddot{\varphi}\dot{\varphi} + 10\dddot{\varphi}\ddot{\varphi} - 10\ddot{\varphi}\dot{\varphi}^3) r\,; & &\quad +(\dddddot{\varphi} - 10\dddot{\varphi}\dot{\varphi}^2 - 15\ddot{\varphi}^2\dot{\varphi} + \dot{\varphi}^5) r \\
&\vdots & &\vdots
\end{aligned}
\qquad (2.3)
$$

## 2.1 Trajectory of a Moving Point

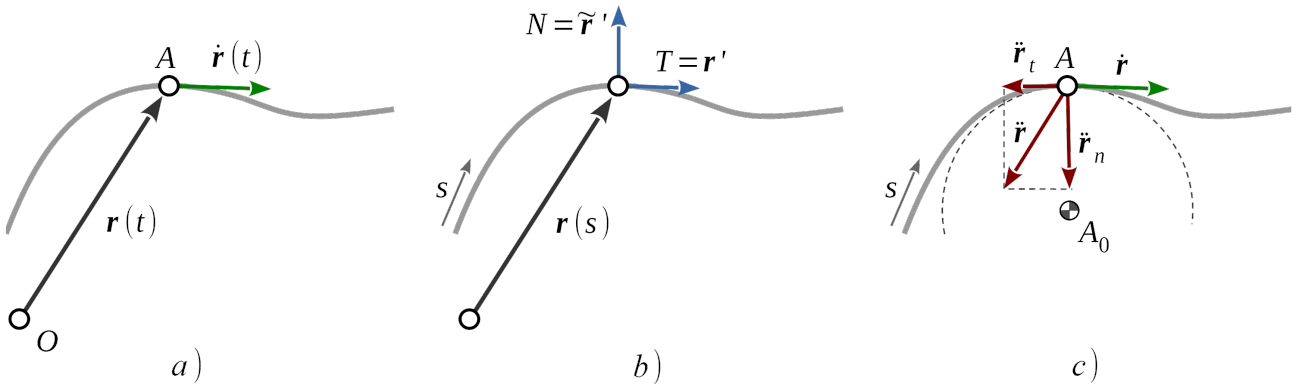

Fig. 2: Trajectory of a moving point

The trajectory of moving point $A$ is described by the vector $\mathbf{r}(t)$ from fixed point $O$ to that point. Its velocity $\dot{\mathbf{r}}(t)$ is always tangential to the trajectory (Fig. 2a). Parametric curve $\mathbf{r}(t)$ is a sufficiently times continuously differentiable function in $\mathbb{R}^2$ [15].

$$\mathbf{r}(t) = \begin{pmatrix} r_1(t) \\ r_2(t) \end{pmatrix}, \quad t \in \mathbb{R}$$

The tangent vector at location $t$ is

$$\dot{\mathbf{r}}(t) = \begin{pmatrix} \dot{r}_1(t) \\ \dot{r}_2(t) \end{pmatrix}$$

> *Remark:*
> The curve parameter $t$ does not necessarily have to be interpreted as time here, but rather as any strictly monotonically increasing curve parameter.

## 2.2 Parameterizing by Arc Length

Introducing *arc length $s$* as curve parameter, an arc element of infinitesimal length $ds$ at location $t$ is given by [2]

$$ds = \sqrt{\dot{\mathbf{r}}^2}\, dt = \sqrt{\dot{r}_1^2 + \dot{r}_2^2}\, dt = \dot{r}\, dt \tag{2.4}$$

Thus we yield the overall arc length $s(t)$ from start location $t_0$ to $t$ by integration

$$s(t) = \int_{t_0}^{t} \sqrt{\dot{\mathbf{r}}^2}\, dt$$

Now we are able to parametrize the curve by arc length, i.e. $\mathbf{r}(s)$

$$\dot{\mathbf{r}} = \frac{d\mathbf{r}}{dt} = \frac{d\mathbf{r}}{ds}\frac{ds}{dt} = \mathbf{r}'\frac{ds}{dt} \tag{2.5}$$

> *Remark:*
> Rewriting this as $\mathbf{r}' = \dot{\mathbf{r}}\dfrac{dt}{ds}$ and using (2.4) proves $\mathbf{r}' = \dfrac{\dot{\mathbf{r}}}{\dot{r}}$ to be a unit vector.

## 2.3 Frenet Formulas in the Plane

We would like to rename the tangential unit vector $\mathbf{r}'$ to $\mathbf{T}$.

> **Lemma 2.2:**
> The orthonormal base $\{\mathbf{T}, \mathbf{N}\}$ is the intrinsic *Frenet coordinate system* for plane curves. The equations
>
> $$\begin{pmatrix} \mathbf{T}' \\ \mathbf{N}' \end{pmatrix} = \kappa \underbrace{\begin{pmatrix} 0 & 1 \\ -1 & 0 \end{pmatrix}}_{\mathbf{J}^T} \cdot \begin{pmatrix} \mathbf{T} \\ \mathbf{N} \end{pmatrix} \tag{2.6}$$
>
> are the *Frenet* equations of the curve $\mathbf{r}(s)$. Term $\kappa$ is called *curvature* [2].

*Proof.* Differentiating $\mathbf{T}^2 = 1$ with respect to curve parameter $s$ gives us

$$\frac{d}{ds}\mathbf{T}^2 = 2\mathbf{T} \cdot \mathbf{T}' = 0.$$

So $\mathbf{T}'$ must be perpendicular to $\mathbf{T}$, thus

$$\mathbf{T}' = \kappa(\mathbf{JT}) = \kappa\tilde{\mathbf{T}} = \kappa\mathbf{N}.$$

with some factor $\kappa$ and normal unit vector $\mathbf{N}$. Further

$$\mathbf{N}' = \mathbf{JT}' = -\kappa\mathbf{T}.$$

So tangent unit vector $\mathbf{T}$ and normal unit vector $\mathbf{N}$ build a right-handed system (Fig.1b). □

> *Remark:*
> The Frenet coordinate system $\{\mathbf{T}, \mathbf{N}\}$ for plane curves conforms to the symplectic orthonormal base $\{\mathbf{r}', \tilde{\mathbf{r}}'\}$.

## 2.4 Center of Curvature

> **Theorem 2.3**:
> The center of curvature $A_0$ of the trajectory at moving point $A$ is given by
>
> $$\mathbf{r}_{AA_0} = \frac{\dot{r}^2}{\tilde{\dot{\mathbf{r}}}\ddot{\mathbf{r}}}\dot{\tilde{\mathbf{r}}} \qquad (2.7)$$
>
> with known velocity $\dot{\mathbf{r}}$ and acceleration $\ddot{\mathbf{r}}$ of point $A$.

*Proof.* We can interpret velocity and normal acceleration of point $A$ as the result of an instantanious rotation about its center of curvature $A_0$ with angular velocity $\dot{\varphi}$ according to equations (2.3) and Figure 1c, i.e.

$$\dot{\mathbf{r}} = \dot{\varphi}\,\tilde{\mathbf{r}}_{A_0 A} \quad and \quad \ddot{\mathbf{r}}_n = -\dot{\varphi}^2\,\mathbf{r}_{A_0 A}. \qquad (2.8)$$

From these two expressions $\ddot{\mathbf{r}}_n = \dot{\varphi}\,\dot{\tilde{\mathbf{r}}}$ can be synthesized and - after multiplying this by $\dot{\tilde{\mathbf{r}}}$ - we get

$$\dot{\varphi} = \frac{\dot{\tilde{\mathbf{r}}}\,\ddot{\mathbf{r}}_n}{\dot{r}^2}.$$

That angular velocity is introduced back into $\dot{\mathbf{r}}$ from (2.8), while being allowed to equate $\dot{\tilde{\mathbf{r}}}\ddot{\mathbf{r}}_n = \dot{\tilde{\mathbf{r}}}\ddot{\mathbf{r}}$ due to the projective character of the scalar product. This finally leads to the location $A_0$ of the center of curvature observed from point $A$ by equation (2.7). □

> *Remark:*
> Equation (2.7) is the vectorial equivalence to the well known scalar formula for the radius of curvature [8,11,15,16]
>
> $$\rho = \frac{(\dot{x}^2 + \dot{y}^2)^{\frac{3}{2}}}{\dot{x}\ddot{y} - \dot{y}\ddot{x}}$$

> **Corrolary 2.4**:
> The magnitude of the curvature vector $\mathbf{r}_{AA_0}$ in equation (2.7) corresponds to the reciprocal of $\kappa$ in Lemma 2.2
>
> $$r_{AA_0} = \frac{1}{\kappa} \qquad (2.9)$$

*Proof.* We substitute $\dot{\mathbf{r}} = \mathbf{r}'\dot{s}$, $\ddot{\mathbf{r}} = \mathbf{r}''\dot{s}^2 + \mathbf{r}'\ddot{s}$ into equation (2.7) according to (2.5) and also use $\mathbf{r}'' = \kappa\tilde{\mathbf{r}}'$ from (2.6).

$$\mathbf{r}_{AA_0} = \frac{\mathbf{r}'^2 \dot{s}^2}{\tilde{\mathbf{r}}'\dot{s}(\kappa\tilde{\mathbf{r}}'\dot{s}^2 + \mathbf{r}'\ddot{s})}\tilde{\mathbf{r}}'\dot{s}$$

Outmultiplying the denominator eliminates the term containing $\ddot{s}$ and taking $\mathbf{r}'^2 = 1$ into account results in $\mathbf{r}_{AA_0} = \frac{\tilde{\mathbf{r}}'}{\kappa}$. Multiplication by $\mathbf{r}'$ yields scalar expression (2.9), which serves as a proof for the term $\kappa$ in expressions (2.6) being the curvature at the point of the the trajectory under consideration. □

# 3. Planar Rigid Body Kinematics

## 3.1 Rigid Body Frames

Traditionally rigid body kinematics uses *Euclidean frames* attached to the bodies. So each body is given a local origin along with an orthonormal basis – $\{x, y\}$ here in $\mathbb{R}^2$. Motion of that frame is then described with respect to a *frame of reference* $\{X, Y\}$ (Fig.3a). The relation between absolute and relative position coordinates reads [3,4,6,12,14,15,17,18]

$$\begin{aligned} X &= a + x\cos\theta - y\sin\theta \\ Y &= b + x\sin\theta + y\cos\theta \end{aligned} \tag{3.1}$$

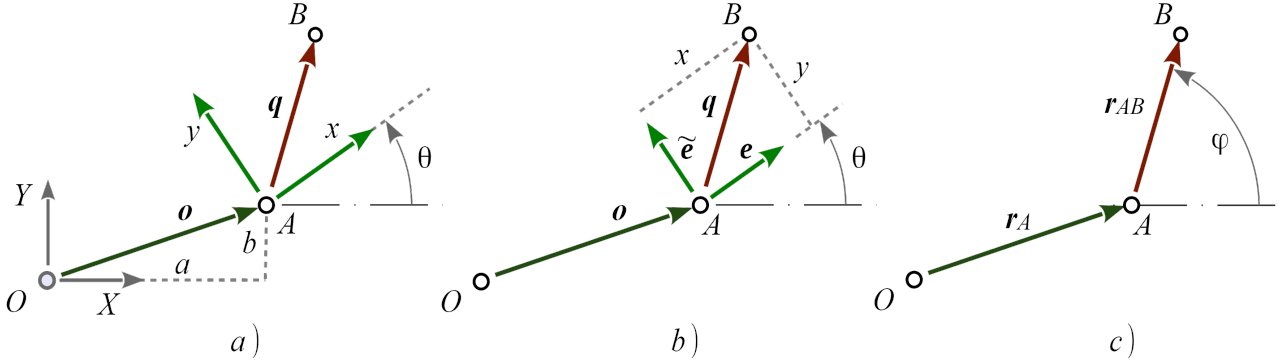

Fig. 3: Moving body frame.

The scalar components in equations (3.1) can be collected in vectors and a skew-symmetric matrix in two different ways.

$$\begin{aligned}
\underbrace{\begin{pmatrix} X \\ Y \end{pmatrix}}_{\mathbf{r}} &= \underbrace{\begin{pmatrix} a \\ b \end{pmatrix}}_{\mathbf{o}} + \underbrace{\begin{pmatrix} \cos\theta & -\sin\theta \\ \sin\theta & \cos\theta \end{pmatrix}}_{\mathbf{R}(\theta)} \cdot \underbrace{\begin{pmatrix} x \\ y \end{pmatrix}}_{\mathbf{q}} \\
&= \underbrace{\begin{pmatrix} a \\ b \end{pmatrix}}_{\mathbf{o}} + \underbrace{\begin{pmatrix} x & -y \\ y & x \end{pmatrix}}_{\mathbf{Q}} \cdot \underbrace{\begin{pmatrix} \cos\theta \\ \sin\theta \end{pmatrix}}_{\mathbf{e}(\theta)}
\end{aligned} \tag{3.2}$$

> *Remark:*
> Orthogonal matrix $\mathbf{R}$ is the *Lie group* of the planar rotation $S\mathbb{O}_2$ [4]. Matrix $\mathbf{R}$ and vectors $\mathbf{o}$ and $\mathbf{e}$ are *kinematic invariants*, as they do not depend on the choice of a body point.

Both matrices $\mathbf{R}$ and $\mathbf{Q}$ can be decomposed into a symmetric and skew-symmetric part, i.e. $\mathbf{R}(\theta) = \cos\theta\,\mathbf{I} + \sin\theta\,\mathbf{J}$ and $\mathbf{Q} = x\,\mathbf{I} + y\,\mathbf{J}$. As already discussed in section 1.2, this gives us a coordinate- and matrix-free notation of equation (3.2) (Fig. 3b)

$$\mathbf{r} = \mathbf{o} + \cos\theta\,\mathbf{q} + \sin\theta\,\tilde{\mathbf{q}} = \mathbf{o} + x\,\mathbf{e} + y\,\tilde{\mathbf{e}}. \tag{3.3}$$

The third and simplest approach is to initially dispense with a frame and start with a vector that is required anyway. This allows the vector $\mathbf{r}_{AB}$ in Figure 3c to be extended to a symplectic orthogonal basis $\{\mathbf{r}_{AB}, \tilde{\mathbf{r}}_{AB}\}$ if necessary. That pragmatic approach is preferred below.

## 3.2 Kinematic Equations

Consider a planar rigid body, which is uniquely defined by two points $A$ and $B$ on it. Position of point $A$ is given with respect to an arbitrary origin $O$, as well as the relative position of point $B$ with respect to $A$. The position of point $B$ is obtained by

$$\mathbf{r}_B = \mathbf{r}_A + \mathbf{r}_{AB}. \tag{3.4}$$

> **Theorem 3.1**:
> The motion of a rigid body in the plane with two points $A$ and $B$ is a combination of the translational motion of point A and the superimposed rotational motion of B around A, which itself is a composition of a *radial* and a *tangential* angular component $\Omega_r$ and $\Omega_t$.
>
> The $k+1$ time derivatives of these last two components are a linear recurrence relation of their previous *radial* and a *tangential* $k$th time derivative components with respect to the symplectic basis $\{\mathbf{r}_{AB}, \tilde{\mathbf{r}}_{AB}\}$ .
>
> Derivation of order $k+1$ of the equation of motion of point $B$ with given motion of point $A$ is
>
> $$\mathbf{r}_B^{(k+1)} = \mathbf{r}_A^{(k+1)} - \underbrace{(\dot{\Omega}_r^{(k)} + \omega\,\Omega_t^{(k)})}_{\Omega_r^{(k+1)}}\cdot \mathbf{r}_{AB} + \underbrace{(\dot{\Omega}_t^{(k)} - \omega\,\Omega_r^{(k)})}_{\Omega_t^{(k+1)}}\cdot \tilde{\mathbf{r}}_{AB} \tag{3.5}$$

*Proof.* The first three time derivatives of equation (3.4) using $\omega = \dot{\varphi}$ are

$$\begin{aligned}
\dot{\mathbf{r}}_B &= \dot{\mathbf{r}}_A + \omega\,\tilde{\mathbf{r}}_{AB} \\
\ddot{\mathbf{r}}_B &= \ddot{\mathbf{r}}_A - \omega^2\,\mathbf{r}_{AB} + \dot{\omega}\,\tilde{\mathbf{r}}_{AB} \\
\dddot{\mathbf{r}}_B &= \dddot{\mathbf{r}}_A - 3\dot{\omega}\omega\,\mathbf{r}_{AB} + (\ddot{\omega} - \omega^3)\,\tilde{\mathbf{r}}_{AB}
\end{aligned} \tag{3.6}$$

We write down a generalised formulation for derivation order $k$.

$$\mathbf{r}_B^{(k)} = \mathbf{r}_A^{(k)} - \Omega_r^{(k)}\,\mathbf{r}_{AB} + \Omega_t^{(k)}\,\tilde{\mathbf{r}}_{AB}$$

Using a negative sign with $\Omega_r^{(k)}$, we take into account the fact that the dominating centripetal acceleration of any order is always directed opposite to vector $\mathbf{r}_{AB}$. Deriving that equation we obtain the accelerations of order $k+1$

$$\mathbf{r}_B^{(k+1)} = \mathbf{r}_A^{(k+1)} - \dot{\Omega}_r^{(k)}\,\mathbf{r}_{AB} - \omega\,\Omega_r^{(k)}\,\tilde{\mathbf{r}}_{AB} + \dot{\Omega}_t^{(k)}\,\tilde{\mathbf{r}}_{AB} - \omega\,\Omega_t^{(k)}\,\mathbf{r}_{AB}$$

which is identical to equation (3.5). □

Calculating the first five time derivatives using recursion (3.5) results in the following radial and tangential angular components:

$$\begin{aligned}
\Omega_r^{(0)} &= -1 \,; & \Omega_t^{(0)} &= 0 \\
\Omega_r^{(1)} &= 0 \,; & \Omega_t^{(1)} &= \omega \\
\Omega_r^{(2)} &= \omega^2 \,; & \Omega_t^{(2)} &= \dot\omega \\
\Omega_r^{(3)} &= 3\dot\omega\omega \,; & \Omega_t^{(3)} &= \ddot\omega - \omega^3 \\
\Omega_r^{(4)} &= 4\ddot\omega\omega + 3\dot\omega^2 - \omega^4 \,; & \Omega_t^{(4)} &= \dddot\omega - 6\dot\omega\omega^2 \\
\Omega_r^{(5)} &= 5\dddot\omega\omega + 10\ddot\omega\dot\omega - 10\dot\omega\omega^3 \,; & \Omega_t^{(5)} &= \ddddot\omega - 10\ddot\omega\omega^2 - 15\dot\omega^2\omega + \omega^5 \\
&\vdots & &\vdots
\end{aligned} \qquad (3.7)$$

> *Remarks:*
>
> 1. The higher-order angular acceleration values $\Omega_r^{(k)}$ and $\Omega_t^{(k)}$ are kinematic invariants. They result from $a_\parallel^{(k)}$ and $a_\perp^{(k)}$ in equations (2.2) and (2.3), respectively, by omitting all summands containing time derivations of magnitude $r$.
> 2. Based on the initial values $\Omega_r^{(0)}$ and $\Omega_t^{(0)}$ and the recursion rule (3.5), the radial angular component $\Omega_r^{(k)}$ is independent of derivative $\omega^{(k-1)}$, just as the tangential angular component $\Omega_t^{(k)}$ is independent of derivative $\omega^{(k-2)}$.

Similar recurrence formulas were presented by Condurache et. al. [4]. The statements in Remark 2 were mentioned by Figliolini et al. in the context of higher-order Bresse circles [7].

# 4. Poles of the Planar Motion

## 4.1 Velocity Pole

> **Theorem 4.1**:
> The general planar motion of a rigid body can be interpreted as an instantaneous pure rotation about a specific point – the *instantaneous center of velocity* $P$. The location of that velocity pole observed from body point $A$ is
>
> $$\mathbf{r}_{AP} = \frac{\dot{\mathbf{r}}_A}{\omega} \qquad (4.1)$$

*Proof.* Substituting point $B$ in the first equation of (3.6) by the pole $P$ having zero velocity, i.e. $\dot{\mathbf{r}}_P = \dot{\mathbf{r}}_A + \omega \cdot \tilde{\mathbf{r}}_{AP} = \mathbf{0}$ resolves to the velocity pole's location (4.1). □

> *Remark:*
> The instantaneous center of velocity is also called *velocity pole* or simply *pole* $P$. Due to its high importance in plane kinematics, we will refer to the velocity pole location as vector $\mathbf{p}$ for the sake of further intensive use.

## 4.2 Pole Acceleration $\ddot{\mathbf{p}}$

While the body point that falls into pole $P$ currently has no velocity, its acceleration is not zero. Pole acceleration is given by the second equation in (3.6).

$$\ddot{\mathbf{p}} = \ddot{\mathbf{r}}_A - \omega^2\,\mathbf{r}_{AP} + \dot\omega\,\tilde{\mathbf{r}}_{AP}\,. \qquad (4.2)$$

## 4.3 Pole Displacement Velocity $\mathbf{u}$

The body point falling into the pole $P$ has zero velocity at current. But the pole – as a virtual point not bound to the body – will change its location over time. The velocity resulting from this displacement is called *pole displacement velocity* [5].

> **Theorem 4.2**:
> In a planar motion the pole displacement velocity $\mathbf{u}$ is directed perpendicular to the pole's acceleration $\ddot{\mathbf{p}}$.
>
> $$\mathbf{u} = \frac{\ddot{\mathbf{p}}}{\omega} \qquad (4.3)$$

*Proof.* We take equation (4.1), apply orthogonal operator $\mathbf{J}$ and substitute vector $\mathbf{r}_{AP} = \mathbf{p} - \mathbf{r}_A$ according to equation (3.4).

$$\dot{\mathbf{r}}_A = -\omega\,\tilde{\mathbf{r}}_{AP} = \omega(\tilde{\mathbf{r}}_A - \tilde{\mathbf{p}})$$

Yet differentiating that equation w.r.t. time yields $\ddot{\mathbf{r}}_A = -\dot{\omega}\,\tilde{\mathbf{r}}_{AP} + \omega\,\dot{\tilde{\mathbf{r}}}_A - \omega\,\dot{\tilde{\mathbf{p}}}$. Reusing $\dot{\tilde{\mathbf{r}}}_A = \omega\,\mathbf{r}_{AP}$ from equation (4.1) again and identifying $\dot{\mathbf{p}}$ as the pole displacement velocity $\mathbf{u}$ results in

$$\ddot{\mathbf{r}}_A = -\dot{\omega}\,\tilde{\mathbf{r}}_{AP} + \omega^2\,\mathbf{r}_{AP} - \omega\,\tilde{\mathbf{u}}\,.$$

Finally comparing the summands herein with those in equation (4.2) leads to the insight, that $\ddot{\mathbf{p}} = -\omega\,\tilde{\mathbf{u}}$. □

> *Remark:*
> The direction of the *pole displacement velocity* $\mathbf{u}$ coincides with the direction of the *pole path tangent* $t$, whereas the *pole acceleration* $\ddot{\mathbf{p}}$ coincides with the direction of the *pole path normal* $n$.

## 4.4 Higher-Order Acceleration Poles

> **Theorem 4.3**:
> In a planar motion there is a specific pole $P_k$ for which the acceleration of order $k$ is zero. That acceleration pole of order $k$ observed from point $A$ is obtained by
>
> $$\mathbf{r}_{AP_k} = \frac{\Omega_r^{(k)} \cdot \mathbf{r}_A^{(k)} + \Omega_t^{(k)} \tilde{\mathbf{r}}_A^{(k)}}{\Omega_r^{(k)2} + \Omega_t^{(k)2}} \qquad (4.4)$$
>
> where the translational acceleration $\mathbf{r}_A^{(k)}$ of some point $A$ as well as the body's radial and tangential angular acceleration state $\Omega_r^{(k)}$ and $\Omega_t^{(k)}$ of order $k$ is known.

*Proof.* We substitute point $B$ in equation (3.4) by acceleration pole $P_k$ and demand it to zero.

$$\mathbf{r}_{P_k}^{(k)} = \mathbf{r}_A^{(k)} - \Omega_r^{(k)} \cdot \mathbf{r}_{AP_k} + \Omega_t^{(k)} \cdot \tilde{\mathbf{r}}_{AP_k} = \mathbf{0}$$

Resolving this equation for $\mathbf{r}_A^{(k)}$, then adding that equation's skew-orthogonal complement while multiplying the first by $\Omega_r^{(k)}$ and the latter by $\Omega_t^{(k)}$

$$\begin{aligned} \mathbf{r}_A^{(k)} &= \Omega_r^{(k)} \cdot \mathbf{r}_{AP_k} - \Omega_t^{(k)} \cdot \tilde{\mathbf{r}}_{AP_k} \mid \cdot \Omega_r^{(k)} \\ + \quad \tilde{\mathbf{r}}_A^{(k)} &= \Omega_r^{(k)} \cdot \tilde{\mathbf{r}}_{AP_k} + \Omega_t^{(k)} \cdot \mathbf{r}_{AP_k} \mid \cdot \Omega_t^{(k)} \end{aligned}$$

results in $\Omega_r^{(k)} \mathbf{r}_A^{(k)} + \Omega_t^{(k)} \tilde{\mathbf{r}}_A^{(k)} = (\Omega_r^{(k)2} + \Omega_t^{(k)2})\mathbf{r}_{AP_k}$. This resolves to the location of the acceleration pole (of order $k$) observed from point $A$ in equation (4.4). □

Inserting the concrete radial and tangential angular acceleration from (3.7) gives

$$\mathbf{r}_{AP_1} = \frac{\dot{\mathbf{r}}_A}{\omega}, \quad \mathbf{r}_{AP_2} = \frac{\omega^2 \ddot{\mathbf{r}}_A + \dot{\omega}\dddot{\mathbf{r}}_A}{\omega^4 + \dot{\omega}^2}, \quad \mathbf{r}_{AP_3} = \frac{3\dot{\omega}\omega\dddot{\mathbf{r}}_A + (\ddot{\omega} - \omega^3)\ddddot{\mathbf{r}}_A}{(3\dot{\omega}\omega)^2 + (\ddot{\omega} - \omega^3)^2},$$

$$\mathbf{r}_{AP_4} = \frac{(4\ddot{\omega}\omega + 3\dot{\omega}^2 - \omega^4)\ddddot{\mathbf{r}}_A + (\dddot{\omega} - 6\dot{\omega}\omega^2)\dddddot{\mathbf{r}}_A}{(4\ddot{\omega}\omega + 3\dot{\omega}^2 - \omega^4)^2 + (\dddot{\omega} - 6\dot{\omega}\omega^2)^2}, \quad \ldots \tag{4.5}$$

> *Remark:*
> The pole equations (4.4) and (4.5) are not kinematic invariants, as they depend on the choice of the body point $A$. This can be changed by simply replacing $A$ with the velocity pole $P$, which we will do below.

## 5. Relative Motion

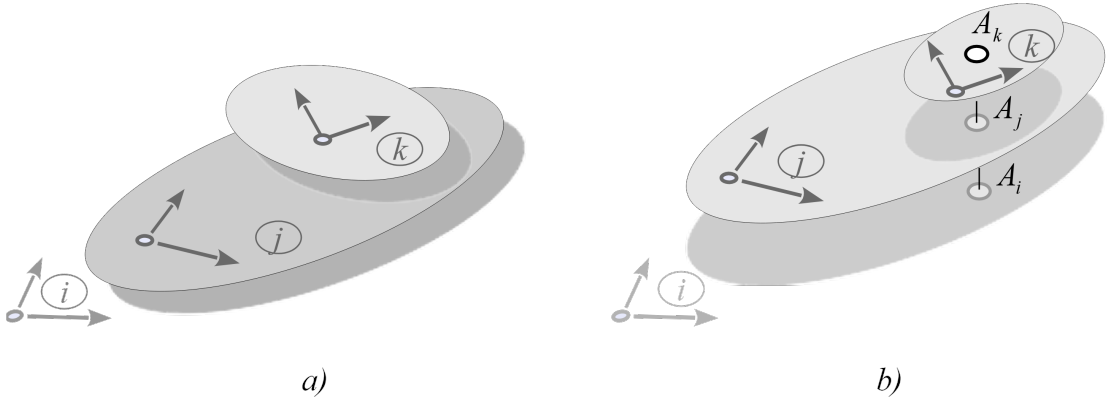

Fig. 4: Relative motion of three bodies.

Let three frames tagged $i$, $j$ and $k$ move relative to each other. We want to consider frame $i$ fixed, which is for better intuition only, without loss of generality though.

### *5.1 Relative Kinematics*

The relativ angle of frame $i$ with respect to frame $j$ is denoted by $\varphi_{ij}$ (Figure 4a). When changing perspective we discover that $\varphi_{ji} = -\varphi_{ij}$.

> **Theorem 5.1**:
> In the case of planar relative motion of three frames $i$, $j$, $k$, their angles, angular velocities and angular accelerations of order $m$ obey the relationship [9]
>
> $$\omega_{ij}^{(m)} + \omega_{jk}^{(m)} + \omega_{ki}^{(m)} = 0 \tag{5.1}$$

*Proof.* Angle $\varphi_{ki}$ in Figure 2a can be expressed via the sum $\varphi_{ki} = \varphi_{ji} + \varphi_{kj}$, which we bring into the form $\varphi_{ij} + \varphi_{jk} + \varphi_{ki} = 0$. If we derive this several times with respect to time, we obtain by substituting $\omega = \dot{\varphi}$

$$\begin{aligned} \varphi_{ij} + \varphi_{jk} + \varphi_{ki} &= 0 \\ \omega_{ij} + \omega_{jk} + \omega_{ki} &= 0 \\ \dot{\omega}_{ij} + \dot{\omega}_{jk} + \dot{\omega}_{ki} &= 0 \\ &\vdots \end{aligned} \tag{5.2}$$

the general equation (5.1) of three frames for their angular accelerations of order $m$. □

> **Remark:**
> Read the term $\dot{\mathbf{r}}_{Aij}$ as "*velocity of point A fixed to frame i with respect to frame j*".

According to Figure 4b the relative velocity $\dot{\mathbf{r}}_{Aki}$ of point A is equal to the sum of the velocities $\dot{\mathbf{r}}_{Akj} + \dot{\mathbf{r}}_{Aji}$. We further observe that $\dot{\mathbf{r}}_{Aij} = -\dot{\mathbf{r}}_{Aji}$. Thus

$$\dot{\mathbf{r}}_{Aij} + \dot{\mathbf{r}}_{Ajk} + \dot{\mathbf{r}}_{Aki} = \mathbf{0}. \tag{5.3}$$

Please take note that this relationship cannot be generalized to accelerations.

> **Remark:**
> Observe the cyclic permutation of indices in equations (5.1), (5.2) and (5.3).

## 5.2 Relative Poles

Equation (4.1) was derived by interpreting velocities with respect to the ground frame. Generalizing from that special case to the relative motion of point $A$ fixed to frame $i$ with respect to frame $j$ we get to the relative velocity pole location $\mathbf{r}_{APij}$ observed from point $A$

$$\mathbf{r}_{APij} = \frac{\tilde{\dot{\mathbf{r}}}_{Aij}}{\omega_{ij}} \tag{5.4}$$

> **Remark:**
> Since $\omega_{ji} = -\omega_{ij}$ and $\dot{\mathbf{r}}_{Aij} = -\dot{\mathbf{r}}_{Aji}$ we verify from expression (5.4) that $\mathbf{r}_{APij} = \mathbf{r}_{APji}$.

With three frames $i, j, k$ we have three relative poles $\mathbf{r}_{APji}$, $\mathbf{r}_{APkj}$ and $\mathbf{r}_{APki}$.

> **Theorem 5.2**: Aronhold-Kennedy
> Any three frames in a planar relative motion have three corresponding relative poles which are collinear.

*Proof.* Collinearity of $\mathbf{r}_{APji}, \mathbf{r}_{APki}, \mathbf{r}_{APkj}$ demands

$$(\mathbf{r}_{APki} - \mathbf{r}_{APji})(\tilde{\mathbf{r}}_{APki} - \tilde{\mathbf{r}}_{APkj}) = 0$$

Reusing equation (5.4) herein leads to

$$\left(\frac{\tilde{\dot{\mathbf{r}}}_{Aki}}{\omega_{ki}} - \frac{\tilde{\dot{\mathbf{r}}}_{Aji}}{\omega_{ji}}\right)\left(\frac{\dot{\mathbf{r}}_{Akj}}{\omega_{kj}} - \frac{\dot{\mathbf{r}}_{Aki}}{\omega_{ki}}\right) = 0.$$

Substituting $\dot{\mathbf{r}}_{Aki} = \dot{\mathbf{r}}_{Aji} + \dot{\mathbf{r}}_{Akj}$ according to (5.3) we get

$$\left(\frac{\tilde{\dot{\mathbf{r}}}_{Aji}}{\omega_{ki}} + \frac{\tilde{\dot{\mathbf{r}}}_{Akj}}{\omega_{ki}} - \frac{\tilde{\dot{\mathbf{r}}}_{Aji}}{\omega_{ji}}\right)\left(\frac{\dot{\mathbf{r}}_{Akj}}{\omega_{kj}} - \frac{\dot{\mathbf{r}}_{Aji}}{\omega_{ki}} - \frac{\dot{\mathbf{r}}_{Akj}}{\omega_{ki}}\right) = 0.$$

Outmultiplying results in three remaining summands

$$\frac{\tilde{\dot{\mathbf{r}}}_{Aji}\,\dot{\mathbf{r}}_{Akj}}{\omega_{ki}\,\omega_{kj}} - \frac{\tilde{\dot{\mathbf{r}}}_{Aji}\,\dot{\mathbf{r}}_{Akj}}{\omega_{ji}\,\omega_{kj}} + \frac{\tilde{\dot{\mathbf{r}}}_{Aji}\,\dot{\mathbf{r}}_{Akj}}{\omega_{ji}\,\omega_{ki}} = 0.$$

From this we can put the common nominator outside the brackets and building the main denominator yields

$$\dot{\tilde{\mathbf{r}}}_{Aji} \dot{\mathbf{r}}_{Akj} \frac{\omega_{ji} - \omega_{ki} + \omega_{kj}}{\omega_{ji}\omega_{ki}\omega_{kj}} = 0.$$

This already completes the proof, since the nominator in that last equation conforms to equation (5.1), which is zero [9]. □

# 6. Bresse Circles

A body point moves along its trajectory, with its velocity always aligned tangentially. This conforms to the Frenet coordinate system, which allows us to decompose the point's acceleration vector of order $k$ into a tangential and a normal component (Fig. 5a).

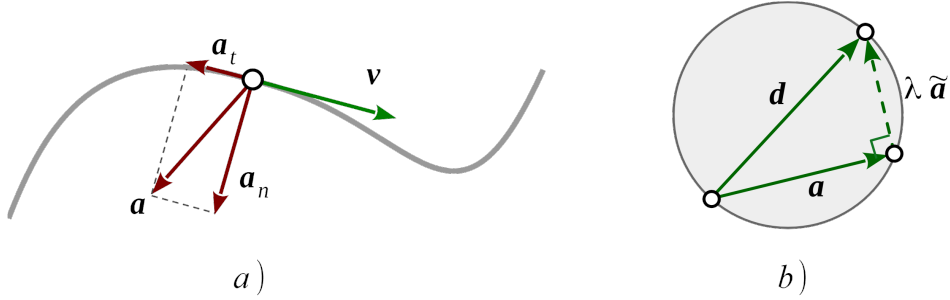

Fig. 5: Tangential and normal acceleration components of a moving body point (a) . Circle equation (b)

From now on, let's denote velocity $\dot{\mathbf{r}}$ by $\mathbf{v}$. The tangential and normal acceleration components of any order $k$ are then

$$\mathbf{r}_T^{(k)} = \frac{\mathbf{v}\mathbf{r}^{(k)}}{v^2}\mathbf{v}, \quad \mathbf{r}_N^{(k)} = \frac{\tilde{\mathbf{v}}\mathbf{r}^{(k)}}{v^2}\tilde{\mathbf{v}}$$

Specifically for points without tangential acceleration, $\mathbf{v}\mathbf{r}^{(k)} = 0$ applies, and correspondingly for points without normal acceleration, $\tilde{\mathbf{v}}\mathbf{r}^{(k)} = 0$ applies.

> **Theorem 6.1**: Bresse Circles of order $k$
>
> The locus of points on a moving plane, that have either no tangential or no normal acceleration of order $k$, is a circle in each case. The respective *zero normal circle* and the *zero tangential circle* are named *first* and *second Bresse circle* of order $k-1$. They are kinematic invariants and their corresponding diameter vectors start from the pole $P$ and point towards the opposite poles $N_{k-1}$ and $T_{k-1}$ respectively.
>
> $$\left. \begin{array}{l} \mathbf{r}_{PN{k-1}} = \dfrac{\mathbf{p}^{(k)}}{\Omega_r^{(k)}} \\[1em] \mathbf{r}_{PT{k-1}} = \dfrac{\tilde{\mathbf{p}}^{(k)}}{\Omega_t^{(k)}} \end{array} \right\} \quad for \ \ k > 1 \qquad (6.1)$$
>
> Both diameter vectors are perpendicular to each other. The acceleration $\mathbf{p}^{(k)}$ of the pole $P$ as well as the body's radial and tangential angular acceleration state $\Omega_r^{(k)}$ and $\Omega_t^{(k)}$ of order $k$ are defined in Theorem 3.1.

*Proof.* We query moving points $Q$ that, observed from pole $P$, satisfy either the condition $\mathbf{v}_Q \cdot \mathbf{r}_Q^{(k)} = 0$ or the condition $\tilde{\mathbf{v}}_Q \cdot \mathbf{r}_Q^{(k)} = 0$. Equation (4.1) gives us $\mathbf{v}_Q = \omega\tilde{\mathbf{r}}_{PQ}$ and equation (3.5) equivalently $\mathbf{r}_Q^{(k)} = \mathbf{p}^{(k)} - \Omega_r^{(k)} \cdot \mathbf{r}_{PQ} + \Omega_t^{(k)} \cdot \tilde{\mathbf{r}}_{PQ}$.

The *zero normal acceleration* condition yields

$$\tilde{\mathbf{v}}_Q\,\mathbf{r}_Q^{(k)} = (-\omega\mathbf{r}_{PQ})(\mathbf{p}^{(k)} - \Omega_r^{(k)}\cdot\mathbf{r}_{PQ} + \Omega_t^{(k)}\cdot\tilde{\mathbf{r}}_{PQ}) = 0 \quad\Rightarrow\quad r_{PQ}^2 - \frac{\mathbf{p}^{(k)}}{\Omega_r^{(k)}}\mathbf{r}_{PQ} = 0 \tag{6.2}$$

and the *zero tangential acceleration* condition results in

$$\tilde{\mathbf{v}}_Q\,\mathbf{r}_Q^{(k)} = (\omega\tilde{\mathbf{r}}_{PQ})(\mathbf{p}^{(k)} - \Omega_r^{(k)}\cdot\mathbf{r}_{PQ} + \Omega_t^{(k)}\cdot\tilde{\mathbf{r}}_{PQ}) = 0 \quad\Rightarrow\quad r_{PQ}^2 - \frac{\tilde{\mathbf{p}}^{(k)}}{\Omega_t^{(k)}}\mathbf{r}_{PQ} = 0 \tag{6.3}$$

From Fig. 5b, we deduce the vector relation $\mathbf{a} + \lambda\tilde{\mathbf{a}} = \mathbf{d}$, where both vectors $\mathbf{a}$ and $\mathbf{d}$ start from the same point on a circle and end at two different points on its circumference, with $\mathbf{d}$ being the diameter vector in particular. Multiplying this equation by $\mathbf{a}$ gives the scalar circle equation

$$\mathbf{a}^2 - \mathbf{d}\mathbf{a} = 0 \tag{6.4}$$

Both of above equations (6.2) and (6.3) have this form of the circle equation. Thus we are able to identify by vector $\mathbf{d}$ the zero normal circle diameter $\mathbf{r}_{PNk-1}$ and the zero tangential circle diameter $\mathbf{r}_{PTk-1}$ in (6.1) from (6.2) and (6.3). The perpendicularity of both to each other is obvious. □

Inserting the concrete acceleration values $\Omega_n^{(k)}$ and $\Omega_t^{(k)}$ from equations (3.7) yields

$$\begin{aligned}
\mathbf{r}_{PN1} &= \frac{\ddot{\mathbf{p}}}{\omega^2}\;;\;\mathbf{r}_{PT1} = \frac{\ddot{\tilde{\mathbf{p}}}}{\dot{\omega}} \\
\mathbf{r}_{PN2} &= \frac{\dddot{\mathbf{p}}}{3\dot{\omega}\omega}\;;\;\mathbf{r}_{PT2} = \frac{\dddot{\tilde{\mathbf{p}}}}{\ddot{\omega} - \omega^3} \\
\mathbf{r}_{PN3} &= \frac{\ddddot{\mathbf{p}}}{4\ddot{\omega}\omega + 3\dot{\omega}^2 - \omega^4}\;;\;\mathbf{r}_{PT3} = \frac{\ddddot{\tilde{\mathbf{p}}}}{\dddot{\omega} - 6\dot{\omega}\omega^2} \\
&\vdots\quad\quad;\quad\quad\vdots
\end{aligned} \tag{6.5}$$

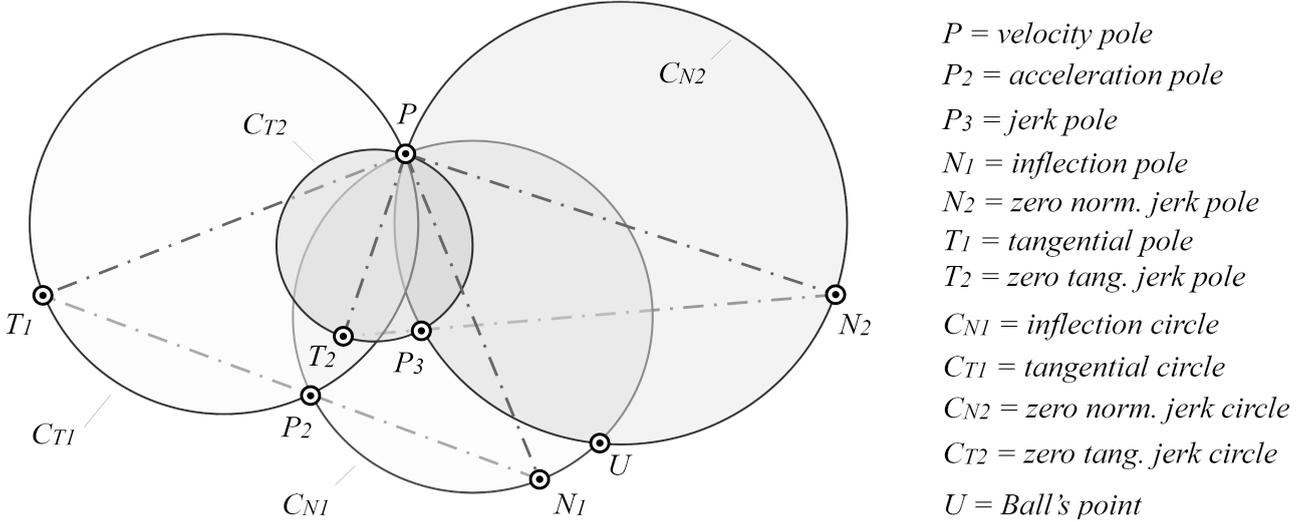

| | |
|---|---|
| $P$ | = velocity pole |
| $P_2$ | = acceleration pole |
| $P_3$ | = jerk pole |
| $N_1$ | = inflection pole |
| $N_2$ | = zero norm. jerk pole |
| $T_1$ | = tangential pole |
| $T_2$ | = zero tang. jerk pole |
| $C_{N1}$ | = inflection circle |
| $C_{T1}$ | = tangential circle |
| $C_{N2}$ | = zero norm. jerk circle |
| $C_{T2}$ | = zero tang. jerk circle |
| $U$ | = Ball's point |

Fig. 6: First and second order Bresse circles

The locations of the *zero normal acceleration poles* and the *zero tangential acceleration poles* are shown in Figure 6 up to the second order for an example case.

> ***Remark:***
> The velocity pole $P$ is a point on all Bresse circles, even though it generally has an acceleration unequal to zero, since it satisfies the conditions (6.2) and (6.3) $\dot{\mathbf{p}}\mathbf{p}^{(k)} = \dot{\tilde{\mathbf{p}}}\mathbf{p}^{(k)} = 0$ due to $\dot{\mathbf{p}} = \mathbf{0}$.

> **Remark:**
> The first Bresse circle of order 1 is called *inflection circle* and the point on it opposite to the pole $P$ is the *inflection pole* $W$. Diameter vector $\mathbf{r}_{PN1}$ is named $\mathbf{r}_{PW}$ from now on. Having introduced this, pole displacement velocity (4.3) can be reformulated to
>
> $$\mathbf{u} = \omega \, \tilde{\mathbf{r}}_{PW} \,. \tag{6.6}$$

## 6.1 Higher-Order Acceleration Poles via Bresse Circles

Pole $P_k$ lies at the intersection point of the Bresse circles $C_{Nk}$ and $C_{Tk}$. The vector from pole $P$ to higher-order pole $P_k$ is orthogonal to the vector from pole $T_{k-1}$ to pole $N_{k-1}$, due to right-angled triangle geometry.

> **Remark:**
> The acceleration pole $P_k$ must lie at the intersection point of the Bresse circles $C_{Nk}$ and $C_{Tk}$ since it has neither normal nor tangential acceleration of order $k$.

The following geometric relationships apply:

$$r_{PP_k} = \lambda(\tilde{\mathbf{r}}_{PNk-1} - \tilde{\mathbf{r}}_{PTk-1}) = \mathbf{r}_{PTk-1} + \mu(\mathbf{r}_{PNk-1} - \mathbf{r}_{PTk-1}) \tag{6.7}$$

This can be resolved for $\lambda$ by multiplication with $(\tilde{\mathbf{r}}_{PNk-1} - \tilde{\mathbf{r}}_{PTk-1})$

$$\lambda = \frac{\mathbf{r}_{PTk-1}(\tilde{\mathbf{r}}_{PNk-1} - \tilde{\mathbf{r}}_{PTk-1})}{(\tilde{\mathbf{r}}_{PNk-1} - \tilde{\mathbf{r}}_{PTk-1})^2} \,.$$

Inserting expressions (6.1) obtains

$$\lambda = \frac{\Omega_r^{(k)} \Omega_t^{(k)}}{\Omega_r^{(k)2} + \Omega_t^{(k)2}}$$

and reusing that in equation (6.7) leads to an invariant equation of higher acceleration poles, equivalent to equation (4.4), where point $A$ is this time replaced by pole $P$.

## 6.2 Ball's Point via Bresse Circles

> **Definition 6.2:**
> A point on the moving plane, which has an instantaneous stationary zero curvature of its path, is a point of *undulation* and is called *Ball's point* after Robert Ball [5,6,9,10,11,15].

The condition for Ball's point is that its velocity $\dot{\mathbf{r}}$ is collinear with both its acceleration vector $\ddot{\mathbf{r}}$ and its jerk vector $\dddot{\mathbf{r}}$, i.e. $\dot{\tilde{\mathbf{r}}}\ddot{\mathbf{r}} = \dot{\tilde{\mathbf{r}}}\dddot{\mathbf{r}} = \mathbf{0}$. Therefore Ball's point $U$ is found to be the other intersection point between the zero-normal circles of first and second order beside the pole.

Using the vector from $N_2$ to $N_1$, Ball's point location observed from the pole is perpendicular to it due to right-angled triangle geometry (Fig. 6) and can be derived in a similar way to equation (6.6).

So from
$$\mathbf{r}_{PU} = \lambda(\tilde{\mathbf{r}}_{PN1} - \tilde{\mathbf{r}}_{PN2}) = \mathbf{r}_{PN1} + \mu(\mathbf{r}_{PN1} - \mathbf{r}_{PN2})$$

we obtain
$$\mathbf{r}_{PU} = \frac{\tilde{\mathbf{r}}_{PN1}\mathbf{r}_{PN2}}{(\mathbf{r}_{PN1} - \mathbf{r}_{PN2})^2}(\tilde{\mathbf{r}}_{PN1} - \tilde{\mathbf{r}}_{PN2}). \tag{6.8}$$

The Bresse vectors used in (6.8) are kinematic invariants. A pure geometric equation for Ball's point location is presented below.

# 7. Geometric Kinematics

From now on, let's denote acceleration $\ddot{\mathbf{r}}$ by $\mathbf{a}$.

## *7.1 Conjugate Points*

While a *body* has an *instantaneous center of velocity* $P$ during motion, the trajectory of some point $A$ fixed to that body has an *instant center of curvature* $A_0$. We already found it via equation (2.7), which we want to rewrite using current notation.

$$\mathbf{r}_{AA_0} = \frac{v_A^2}{\tilde{\mathbf{v}}_A \mathbf{a}_A}\tilde{\mathbf{v}}_A \tag{7.1}$$

> **Definition 7.1:**
> A point $A$ on a moving plane and the instant center of curvature $A_0$ of its trajectory with respect to some other plane are called *conjugate points*. Together with the pole $P$ they lie on a common straight line - the *pole ray* [9,10,15].

The fact that both vectors $\mathbf{r}_{AP}$ in equation (4.1) and $\mathbf{r}_{AA_0}$ in equation (7.1) are directed orthogonal to $\mathbf{v}_A$ from point $A$ suffice as a proof of the collinearity of $P$, $A$ and $A_0$.

If you happen to know two pairs of conjugate points, say $\{A, A_0\}$ and $\{B, B_0\}$, you can find the pole in the intersection of both pole rays, i.e. $\mathbf{r}_P = \mathbf{r}_A + \lambda_A \mathbf{r}_{AA_0}$ and $\mathbf{r}_P = \mathbf{r}_B + \lambda_B \mathbf{r}_{BB_0}$. Equating both, removing $\lambda_B$ via multiplication by $\tilde{\mathbf{r}}_{BB_0}$ or removing $\lambda_A$ via multiplication by $\tilde{\mathbf{r}}_{AA_0}$ instead, leads to [9]

$$\mathbf{r}_{AP} = \frac{\tilde{\mathbf{r}}_{BB_0}\mathbf{r}_{AB}}{\tilde{\mathbf{r}}_{AA_0}\mathbf{r}_{BB_0}}\mathbf{r}_{AA_0} \quad and \quad \mathbf{r}_{BP} = \frac{\tilde{\mathbf{r}}_{AA_0}\mathbf{r}_{AB}}{\tilde{\mathbf{r}}_{AA_0}\mathbf{r}_{BB_0}}\mathbf{r}_{BB_0}. \tag{7.2}$$

Equations (7.2) represent a pure geometric way to find the velocity pole location via two pairs of conjugate points.

## *7.2 Euler-Savary Equation*

The famous equation of Euler-Savary is usually found as a scalar equation in the literature. A new, more general vectorial form has been presented by the author of this publication in [9,10,11].

**Theorem 7.2: Euler-Savary**

The equation of Euler-Savary characterizes the curvature of the trajectory of a point $A$ on a moving plane. Its vectorial form reads

$$\mathbf{r}_{AA_0} = \frac{\mathbf{r}_{PA}^2}{\mathbf{r}_{PW}\mathbf{r}_{PA} - \mathbf{r}_{PA}^2}\mathbf{r}_{PA} \tag{7.3}$$

*Proof.* In equation (7.1) of the instant center of curvature at point $A$ we substitute velocity $\mathbf{v}_A = \omega\tilde{\mathbf{r}}_{PA}$ from equation (4.1) and acceleration $\mathbf{a}_A = \mathbf{a}_P - \omega^2\mathbf{r}_{PA} + \dot{\omega}\tilde{\mathbf{r}}_{PA}$ from equation (4.2), resulting in

$$\mathbf{r}_{AA_0} = \frac{\omega^2\mathbf{r}_{PA}}{(-\omega\mathbf{r}_{PA})(\mathbf{a}_P - \omega^2\mathbf{r}_{PA} + \dot{\omega}\tilde{\mathbf{r}}_{PA})}(-\omega\mathbf{r}_{PA})$$

Outmultiplying the denominator, replacing the pole acceleration $\mathbf{a}_P = \omega^2\mathbf{r}_{PW}$ from equation (6.5) and shortening the fraction, we obtain the purely geometric vector form of the Euler-Savary equation (7.3). □

*Remark:*

The well known equivalent *scalar Euler-Savary* equation can be directly derived from equation (7.3) [11]

$$\rho = \frac{r^2}{D\sin\theta - r}. \tag{7.4}$$

It is called the *second form of Euler-Savary equation* by Pennestri [15]. For the meaning of the variables herein also see [11]. In this scalar form, the sign of distance $r$ must be taken into account separately. The vector form (7.3) already contains the necessary directional information.

## 7.3 Geometric Equation of the Inflection Pole

If we again happen to know two pairs of conjugate points $\{A, A_0\}$ and $\{B, B_0\}$, we can write down Euler-Savary (7.3) two times

$$\mathbf{r}_{AA_0} = \frac{\mathbf{r}_{PA}^2}{\mathbf{r}_{PW}\mathbf{r}_{PA} - \mathbf{r}_{PA}^2}\mathbf{r}_{PA}; \quad \mathbf{r}_{BB_0} = \frac{\mathbf{r}_{PB}^2}{\mathbf{r}_{PW}\mathbf{r}_{PB} - \mathbf{r}_{PB}^2}\mathbf{r}_{PB}$$

From these two equations the inflection pole location $\mathbf{r}_{PW}$ can be synthesized.

$$\mathbf{r}_{PW} = \frac{1}{\tilde{\mathbf{r}}_{PA}\mathbf{r}_{PB}}\left(r_{PB}^2\left(\frac{r_{PB}^2}{\mathbf{r}_{BB_0}\mathbf{r}_{PB}} + 1\right)\tilde{\mathbf{r}}_{PA} - r_{PA}^2\left(\frac{r_{PA}^2}{\mathbf{r}_{AA_0}\mathbf{r}_{PA}} + 1\right)\tilde{\mathbf{r}}_{PB}\right) \tag{7.5}$$

Equation (7.5) is a new vector equivalent to the construction of Bobillier to find the inflection pole $W$ [11].

## 7.4 Geometric Equation of Ball's Point

A pure geometric vector equation of Ball's point location was found and presented by the author of this paper in [11]. It has been derived via vectorization of Bereis' construction. If we know the inflection pole $W$ beside the pole $P$ and two pairs of conjugate points $\{A, A_0\}$ and $\{B, B_0\}$, we get Ball's point location via

$$\mathbf{r}_{PU} = \frac{\mathbf{r}_{PW}\,\mathbf{r}_{PH}}{r_{PH}^2}\mathbf{r}_{PH} \quad \text{with} \quad \mathbf{r}_{PH} = \frac{(\tilde{\mathbf{r}}_{PW}\,\mathbf{r}_{PA})\mathbf{r}_{BB_0} - (\tilde{\mathbf{r}}_{PW}\,\mathbf{r}_{PB})\mathbf{r}_{AA_0}}{\tilde{\mathbf{r}}_{AA_0}\,\mathbf{r}_{BB_0}}. \tag{7.6}$$

Helper point $H$ and Ball's point $U$ lie on a common pole ray [10,11].

## 7.5 Polodes

The locus of all pole positions in the moving plane is the *moving polode* and the locus of all pole positions in the fixed plane is the *fixed polode* [5,12,14,15,17,18]. To handle these curves in relation to different coordinate systems (Fig. 3a), we use equation (3.2) in the rotation matrix form for the pole vector $\mathbf{p}$ in the fixed plane and pole location vector $\mathbf{q}$ in the moving plane.

$$\mathbf{p} = \mathbf{o} + \mathbf{R}\mathbf{q}. \tag{7.7}$$

Under the pragmatic assumption that the angular velocity of the moving plane is constant $\dot{\theta} = \frac{d\theta}{dt} = 1$, time $t$ can be replaced by angular position $\theta$, allowing us to focus exclusively on the geometric aspects of the motion [3,12,14,15].

> **Theorem 7.3:**
> The equations of the fixed and the moving polode – the latter already rotated into the fixed plane – are
>
> $$\begin{aligned} \mathbf{p} &= \mathbf{o} + \tilde{\mathbf{o}}' \\ \mathbf{R}\mathbf{q} &= \tilde{\mathbf{o}}' \end{aligned} \tag{7.8}$$

*Proof.* Differentiating equation (7.7) with respect to $\theta$ (denoted by prime) yields $\mathbf{p}' = \mathbf{o}' + \mathbf{R}'\mathbf{q}$. For the velocity pole $\mathbf{p}' = \mathbf{0}$ applies. Using the derivative $\mathbf{R}' = \mathbf{R}\mathbf{J}$ according to (1.2) we obtain $\mathbf{R}\mathbf{J}\mathbf{q} = \mathbf{R}\tilde{\mathbf{q}} = -\mathbf{o}'$. Applying the orthogonal operator $-\mathbf{J}$ to that result yields $\mathbf{R}\mathbf{q} = \tilde{\mathbf{o}}'$ and substituting this into equation (7.7) gives us the polode equations (7.8). □

> **Theorem 7.4:**
> In the pole $P$ both polodes have collinear tangents and rolling contact without sliding.

*Proof.* Differentiating both equations (7.8) for $\theta$ gives us the tangents, now to be interpreted as the geometric pole displacement velocity. The first equation in (7.8) yields then

$$\mathbf{p}' = \mathbf{o}' + \tilde{\mathbf{o}}''.$$

Deriving the second equation in (7.8) gives

$$\mathbf{R}\mathbf{q}' + \mathbf{R}'\mathbf{q} = \tilde{\mathbf{o}}''.$$

Applying $\mathbf{R}' = \mathbf{R}\mathbf{J}$ to the second term again obtains $\mathbf{R}'\mathbf{q} = \mathbf{R}\tilde{\mathbf{q}} = -\mathbf{o}'$. So we get for the polode tangents in the pole,

$$\begin{aligned} \mathbf{p}' &= \mathbf{o}' + \tilde{\mathbf{o}}'' \\ \mathbf{R}\mathbf{q}' &= \mathbf{o}' + \tilde{\mathbf{o}}'' = \mathbf{p}' \end{aligned} \tag{7.9}$$

which completes the proof for equality of magnitude and direction of the pole velocity $\mathbf{R}\mathbf{q}' = \mathbf{p}'$ in the fixed plane along both polodes.

The differential arc length $ds = \sqrt{x'(\theta)^2 + y'(\theta)^2}$ of a curve equals the magnitude of its corresponding tangent. The moving and fixed polode do have equal magnitudes of their tangents $\mathbf{p}'$ and $\mathbf{q}'$, which proofs the pure rolling contact of both polodes. □

> **Remark:**
> The first derivative $\mathbf{p}'$ is the geometric pole displacement velocity corresponding to velocity $\mathbf{u}$ from equation (4.3). Its magnitude is the diameter $r_{PW}$ of the inflection circle.

E.A. Dijksman has discussed polodes in detail using complex numbers in [5].

## 7.6 Higher-Order Derivatives of the Polodes

**Theorem 7.5:**
The $k$th derivative of the polodes conforming to higher-order geometric accelerations of the pole $P$ are:

$$\mathbf{p}^{(k)} = \mathbf{o}^{(k)} + \tilde{\mathbf{o}}^{(k+1)}$$
$$\mathbf{R}\mathbf{q}^{(k)} = \sum_{i=0}^{k} a_{i,k}(-\mathbf{J})^i \mathbf{p}^{(k-i)} \quad with \quad a_{i,k} = \frac{k!}{i!(k-i)!} \in \mathbb{N} \tag{7.10}$$

*Proof:* The $k$-th derivative of the fixed polode $\mathbf{p}^{(k)}$ can be predicted from the first equation in (7.9). To treat the moving polode of the second equation in (7.9), we make the following approach

$$\mathbf{R}\mathbf{q}^{(1)} = \mathbf{Q}_1 \quad with \quad \mathbf{Q}_1 = \mathbf{p}^{(1)}$$

Deriving this with respect to $\theta$ is

$$\begin{aligned}\mathbf{R}\mathbf{q}^{(2)} &= \mathbf{Q}_1' - \mathbf{R}'\mathbf{q}^{(1)} \\ &= \mathbf{Q}_1' + (-\mathbf{J})\mathbf{R}\mathbf{q}^{(1)} \\ &= \mathbf{Q}_1' + (-\mathbf{J})\mathbf{Q}_1 \\ &= \underbrace{\mathbf{p}^{(2)} + (-\mathbf{J})\mathbf{p}^{(1)}}_{\mathbf{Q}_2}\end{aligned}$$

The next derivation step yields

$$\begin{aligned}\mathbf{R}\mathbf{q}^{(3)} &= \mathbf{Q}_2' + (-\mathbf{J})\mathbf{Q}_2 \\ &= \underbrace{\mathbf{p}^{(3)} + 2(-\mathbf{J})\mathbf{p}^{(2)} + (-\mathbf{J})^2\mathbf{p}^{(1)}}_{\mathbf{Q}_3}\end{aligned}$$

The $k-1$ derivative obtains

$$\begin{aligned}\mathbf{R}\mathbf{q}^{(k-1)} &= \mathbf{Q}_{k-2}' + (-\mathbf{J})\mathbf{Q}_{k-2} \\ &= \underbrace{\mathbf{p}^{(k-1)} + a_{2,k-1}(-\mathbf{J})\mathbf{p}^{(k-2)} + a_{3,k-1}(-\mathbf{J})^2\mathbf{p}^{(k-3)} + \cdots + a_{k-2,k-1}(-\mathbf{J})^{k-2}\mathbf{p}^{(2)} + (-\mathbf{J})^{k-1}\mathbf{p}^{(1)}}_{\mathbf{Q}_{k-1}}\end{aligned}$$

Finally we get for the $k$th derivation while reusing $\mathbf{Q}_{k-1}$

$$\begin{aligned}\mathbf{R}\mathbf{q}^{(k)} &= \mathbf{Q}_{k-1}' + (-\mathbf{J})\mathbf{Q}_{k-1} \\ &= \mathbf{p}^{(k)} + a_{2,k-1}(-\mathbf{J})\mathbf{p}^{(k-1)} + a_{3,k-1}(-\mathbf{J})^2\mathbf{p}^{(k-2)} + \cdots + (-\mathbf{J})^{k-1}\mathbf{p}^{(2)} + \\ &\quad (-\mathbf{J})\mathbf{p}^{(k-1)} + a_{2,k-1}(-\mathbf{J})^2\mathbf{p}^{(k-2)} + \cdots + a_{k-2,k-1}(-\mathbf{J})^{k-1}\mathbf{p}^{(2)} + (-\mathbf{J})^k\mathbf{p}^{(1)} \\ &= \underbrace{\mathbf{p}^{(k)} + (1 + a_{2,k-1})(-\mathbf{J})\mathbf{p}^{(k-1)} + (a_{2,k-1} + a_{3,k-1})(-\mathbf{J})^2\mathbf{p}^{(k-2)} + \cdots + (a_{k-2,k-1} + 1)(-\mathbf{J})^{k-1}\mathbf{p}^{(2)} + (-\mathbf{J})^k\mathbf{p}^{(1)}}_{\mathbf{Q}_k}\end{aligned}$$

The coefficients $a_{i,k}$ are positive integer values that result from the sum of the preceding, adjacent values $a_{i,k} = a_{i-1,k-1} + a_{i,k-1}$. These are entries of Pascal's triangle, which can be calculated according to the rule $a_{i,k} = \frac{k!}{i!(k-i)!}$ [13]. □

Applying equation (7.10) to the moving polode for several $k$'s with the powers of $-\mathbf{J}$ according to equation (1.2) yields

$$\begin{aligned}
\mathbf{R}\mathbf{q}' &= \mathbf{p}' \\
\mathbf{R}\mathbf{q}'' &= \mathbf{p}'' - \tilde{\mathbf{p}}' \\
\mathbf{R}\mathbf{q}''' &= \mathbf{p}''' - 2\tilde{\mathbf{p}}'' - \mathbf{p}' \\
\mathbf{R}\mathbf{q}^{IV} &= \mathbf{p}^{IV} - 3\tilde{\mathbf{p}}''' - 3\mathbf{p}'' + \tilde{\mathbf{p}}' \\
\mathbf{R}\mathbf{q}^{V} &= \mathbf{p}^{V} - 4\tilde{\mathbf{p}}^{IV} - 6\mathbf{p}''' + 4\tilde{\mathbf{p}}'' + \mathbf{p}' \\
\mathbf{R}\mathbf{q}^{VI} &= \mathbf{p}^{VI} - 5\tilde{\mathbf{p}}^{V} - 10\mathbf{p}^{IV} + 10\tilde{\mathbf{p}}''' + 5\mathbf{p}'' - \tilde{\mathbf{p}}' \\
&\vdots
\end{aligned} \qquad (7.11)$$

## 7.6 Curvature Radii of the Polodes

**Theorem 7.7:**
The equations of the curvature radius vectors of the fixed and moving polodes in their contact point $P$ with respect to the fixed frame are

$$\boldsymbol{\rho}_f = \frac{\mathbf{p}'^2}{\tilde{\mathbf{p}}'\mathbf{p}''}\tilde{\mathbf{p}}'$$

$$\mathbf{R}\boldsymbol{\rho}_m = \frac{\mathbf{p}'^2}{\tilde{\mathbf{p}}'\mathbf{p}'' - \mathbf{p}'^2}\tilde{\mathbf{p}}'$$

(7.12)

*Proof.* The curvature radius of the fixed polode $\boldsymbol{\rho}_f$ in equation (7.12) is obtained by substituting geometric velocity $\mathbf{p}'$ and geometric acceleration $\mathbf{p}''$ into equation (2.7). Similarly, we obtain the radius vector of curvature of the moving polode $\mathbf{R}\boldsymbol{\rho}_m$ in the fixed plane by substituting $\mathbf{R}\mathbf{q}'$ and $\mathbf{R}\mathbf{q}''$ instead.

$$\mathbf{R}\boldsymbol{\rho}_m = \frac{(\mathbf{R}\mathbf{q}')^2}{(\mathbf{R}\tilde{\mathbf{q}}')(\mathbf{R}\mathbf{q}'')}\mathbf{R}\tilde{\mathbf{q}}'$$

Inserting the first two equations from (7.11) herein completes equations (7.12). □

## 7.7 Bottema Invariants

Now we align the fixed and moving origin by setting $\mathbf{o} = \mathbf{0}$ and align their axes by setting $\mathbf{R}(\theta = 0) = \mathbf{I}$. We place the common origin at the pole $P$ with the $X$ axis directed along the common tangent of the polodes in $P$. The fixed and moving frame are then called *canonical systems* [3,6,14,15,17,18].

For the derivatives of $\mathbf{o}$ and its components $a, b$ we adopt the convention of notation from [3], i.e. $@_k = \frac{d^k@}{d\theta^k}$. Therefore, the following applies.

$$\mathbf{o}_0 = \begin{pmatrix} 0 \\ 0 \end{pmatrix}; \quad \mathbf{o}_1 = \begin{pmatrix} 0 \\ 0 \end{pmatrix}; \quad \mathbf{o}_2 = \begin{pmatrix} 0 \\ b_2 \end{pmatrix}; \quad \mathbf{o}_3 = \begin{pmatrix} a_3 \\ b_3 \end{pmatrix}; \quad \cdots \quad \mathbf{o}_k = \begin{pmatrix} a_k \\ b_k \end{pmatrix}$$

*Remark:*
The scalar derivatives $b_2, a_3, b_3, ..., a_k, b_k$ are referred to as *Bottema invariants* after *O. Bottema*, who introduced them for values that are independent of specific body points in [3]. We have used the term *invariants* in a less strict sense.

Pennestri discusses Bottema's instantaneous invariants in great detail [14,15] and outlines a valuable historical context relating to Müller, Cesaro, Krause [12], Veldkamp [17,18], Freudenstein and others.

The equation of the fixed polode and its derivatives up to the second order according to (7.10), (7.11) are

$$\begin{aligned} \mathbf{p} &= \mathbf{o}_0 + \tilde{\mathbf{o}}_1 = \begin{pmatrix} 0 \\ 0 \end{pmatrix} \\ \mathbf{p}' &= \mathbf{o}_1 + \tilde{\mathbf{o}}_2 = \begin{pmatrix} -b_2 \\ 0 \end{pmatrix} \\ \mathbf{p}'' &= \mathbf{o}_2 + \tilde{\mathbf{o}}_3 = \begin{pmatrix} -b_3 \\ b_2 + a_3 \end{pmatrix} \end{aligned} \qquad (7.13)$$

For the moving polode and its derivatives up to the second order, we obtain simultaneously

$$\begin{aligned} \mathbf{q}_p &= \mathbf{p} = \begin{pmatrix} 0 \\ 0 \end{pmatrix} \\ \mathbf{q}'_p &= \mathbf{p}' = \begin{pmatrix} -b_2 \\ 0 \end{pmatrix} \\ \mathbf{q}''_p &= \mathbf{p}'' - \tilde{\mathbf{p}}' = \begin{pmatrix} -b_3 \\ 2b_2 + a_3 \end{pmatrix} \end{aligned} \qquad (7.14)$$

> *Remark:*
> Here $b_2$ is the diameter of the inflection circle.

We obtain the curvature radii vectors of the polodes by inserting derivatives (7.13) into (7.12).

$$\begin{aligned} \boldsymbol{\rho}_f &= \frac{b_2}{a_3 + b_2} \begin{pmatrix} 0 \\ b_2 \end{pmatrix} \\ \mathbf{R}\boldsymbol{\rho}_m &= \frac{b_2}{a_3 + 2b_2} \begin{pmatrix} 0 \\ b_2 \end{pmatrix} \end{aligned} \qquad (7.15)$$

Taking the inverse magnitudes of radius vectors (7.15) gives their curvatures. The difference between these leads to

$$\frac{1}{\rho_m} - \frac{1}{\rho_f} = \frac{a_3 + 2b_2}{b_2^2} - \frac{a_3 + b_2}{b_2^2} = \frac{1}{b_2} \qquad (7.16)$$

the *second Euler-Savary* equation relating the curvature of the polodes and the inflection circle diameter.

# 5. Conclusion

It has been shown that symplectic geometry is excellently suited to kinematic analysis in $\mathbb{R}^2$. The resulting equations are straightforward and mostly intuitive. They do not involve coordinates or matrices, unless we explicitly want them to. The equations contain repeating patterns that are easy to spot:

1. Areal expressions such as $\mathbf{ab}$, $\tilde{\mathbf{a}}\mathbf{b}$, $\mathbf{a}^2$, which often occur in fractions [eqn. (2.7),(6.4),(6.8),(7.1),(7.2), (7.3),(7.5),(7.6),(7.12)].
2. Similarity transformations such as $a\mathbf{c} + b\tilde{\mathbf{c}}$, as a combination of scaling and rotation [eqn. (1.1),(2.2), (3.3),(3.5),(4.2),(4.4),(4..5),(7.5)].

The significance of the general radial and tangential expressions of higher-order angular accelerations $\Omega_r^{(k)}$ and $\Omega_t^{(k)}$ in kinematic equations for rigid bodies is noteworthy. Several proofs and the higher-order equations in plane kinematics are presented here for the first time in this general symplectic form.